\documentstyle[titlepage,12pt]{article}

\newtheorem{thm}{Theorem}[section]
\newtheorem{prop}[thm]{Proposition}
\newtheorem{cor}[thm]{Corollary}
\newtheorem{lem}[thm]{Lemma}
\newtheorem{conj}[thm]{Conjecture}
\newtheorem{exa}[thm]{Example}
\newtheorem{defn}[thm]{Definition}
\newcommand{\ben}{\begin{enumerate}}
\newcommand{\een}{\end{enumerate}}
\newcommand{\ble}{\begin{lem}}
\newcommand{\ele}{\end{lem}}
\newcommand{\bth}{\begin{thm}}
\newcommand{\eth}{\end{thm}}
\newcommand{\bpr}{\begin{prop}}
\newcommand{\epr}{\end{prop}}
\newcommand{\bco}{\begin{cor}}
\newcommand{\eco}{\end{cor}}
\newcommand{\bcon}{\begin{conj}}
\newcommand{\econ}{\end{conj}}
\newcommand{\bde}{\begin{defn}}
\newcommand{\ede}{\end{defn}}
\newcommand{\bex}{\begin{exa}}
\newcommand{\eex}{\end{exa}}
\newcommand{\barr}{\begin{array}}
\newcommand{\earr}{\end{array}}
\newcommand{\btab}{\begin{tabular}}
\newcommand{\etab}{\end{tabular}}
\newcommand{\beq}{\begin{equation}}
\newcommand{\eeq}{\end{equation}}
\newcommand{\bea}{\begin{eqnarray*}}
\newcommand{\eea}{\end{eqnarray*}}
\newcommand{\bce}{\begin{center}}
\newcommand{\ece}{\end{center}}
\newcommand{\bpi}{\begin{picture}}
\newcommand{\epi}{\end{picture}}
\newcommand{\bfi}{\begin{figure} \begin{center}}
\newcommand{\efi}{\end{center} \end{figure}}
\newcommand{\capt}{\caption}
\newcommand{\bsl}{\begin{slide}{}}
\newcommand{\esl}{\end{slide}}

\newcommand{\bib}{thebibliography}
\newcommand{\pf}{{\bf Proof.}}
\newcommand{\qed}{\rule{1ex}{1ex}}
\newcommand{\Qed}{\rule{1ex}{1ex} \medskip}

\newcommand{\ol}{\overline}
\newcommand{\hs}[1]{\hspace{#1}}

\newcommand{\st}[1]{\rule{#1}{0pt}}

\newcommand{\mbc}[1]{\makebox(0,0){#1}}

\newcommand{\jn}{\vee}

\newcommand{\case}[4]{\left\{\barr{ll}#1&\mbox{#2}\\#3&\mbox{#4}\earr\right.}
\newcommand{\flf}[2]{\left\lfloor\frac{#1}{#2}\right\rfloor}

\newcommand{\ra}{\rightarrow}

\newcommand{\la}{\lambda}

\newcommand{\De}{\Delta}

\newcommand{\Wb}{\ol{W}}
\newcommand{\Kb}{\ol{K}}
\newcommand{\kb}{\ol{k}}

\newcommand{\lh}{\hat{l}}

\newcommand{\scz}{\scriptsize}

\newcommand{\BTAa}{\put(0,80){\circle*{3}}}

\newcommand{\BTBa}{\put(-80,60){\circle*{3}}}
\newcommand{\BTBb}{\put(80,60){\circle*{3}}}

\newcommand{\BTDa}{\put(-140,20){\circle*{3}}}
\newcommand{\BTDb}{\put(-100,20){\circle*{3}}}

\newcommand{\BTDd}{\put(-20,20){\circle*{3}}}
\newcommand{\BTDe}{\put(20,20){\circle*{3}}}

\newcommand{\BTDh}{\put(140,20){\circle*{3}}}

\newcommand{\BTEa}{\put(-150,0){\circle*{3}}}
\newcommand{\BTEb}{\put(-130,0){\circle*{3}}}
\newcommand{\BTEc}{\put(-110,0){\circle*{3}}}
\newcommand{\BTEd}{\put(-90,0){\circle*{3}}}

\newcommand{\BTEg}{\put(-30,0){\circle*{3}}}
\newcommand{\BTEh}{\put(-10,0){\circle*{3}}}
\newcommand{\BTEi}{\put(10,0){\circle*{3}}}

\newcommand{\TTAa}{\put(0,90){\circle*{3}}}

\newcommand{\TTBa}{\put(-90,60){\circle*{3}}}
\newcommand{\TTBb}{\put(0,60){\circle*{3}}}
\newcommand{\TTBc}{\put(90,60){\circle*{3}}}

\newcommand{\TTCa}{\put(-120,30){\circle*{3}}}

\newcommand{\TTCc}{\put(-60,30){\circle*{3}}}
\newcommand{\TTCd}{\put(-30,30){\circle*{3}}}

\newcommand{\TTCf}{\put(30,30){\circle*{3}}}
\newcommand{\TTCg}{\put(60,30){\circle*{3}}}

\newcommand{\TTCi}{\put(120,30){\circle*{3}}}

\newcommand{\TTDa}{\put(-130,0){\circle*{3}}}

\newcommand{\TTDc}{\put(-110,0){\circle*{3}}}

\newcommand{\TTDg}{\put(-70,0){\circle*{3}}}

\newcommand{\TTDi}{\put(-50,0){\circle*{3}}}
\newcommand{\TTDj}{\put(-40,0){\circle*{3}}}

\newcommand{\TTDl}{\put(-20,0){\circle*{3}}}

\newcommand{\TTDp}{\put(20,0){\circle*{3}}}

\newcommand{\BTAaa}{\put(0,80){\line(-4,-1){80}}}
\newcommand{\BTAab}{\put(0,80){\line(4,-1){80}}}

\newcommand{\BTBaa}{\put(-80,60){\line(-2,-1){40}}}
\newcommand{\BTBab}{\put(-80,60){\line(2,-1){40}}}
\newcommand{\BTBbc}{\put(80,60){\line(-2,-1){40}}}
\newcommand{\BTBbd}{\put(80,60){\line(2,-1){40}}}

\newcommand{\BTCaa}{\put(-120,40){\line(-1,-1){20}}}

\newcommand{\BTCbd}{\put(-40,40){\line(1,-1){20}}}
\newcommand{\BTCce}{\put(40,40){\line(-1,-1){20}}}

\newcommand{\BTCdh}{\put(120,40){\line(1,-1){20}}}

\newcommand{\BTDaa}{\put(-140,20){\line(-1,-2){10}}}
\newcommand{\BTDab}{\put(-140,20){\line(1,-2){10}}}
\newcommand{\BTDbc}{\put(-100,20){\line(-1,-2){10}}}
\newcommand{\BTDbd}{\put(-100,20){\line(1,-2){10}}}

\newcommand{\BTDdg}{\put(-20,20){\line(-1,-2){10}}}
\newcommand{\BTDdh}{\put(-20,20){\line(1,-2){10}}}
\newcommand{\BTDei}{\put(20,20){\line(-1,-2){10}}}

\newcommand{\TTAaa}{\put(0,90){\line(-3,-1){90}}}
\newcommand{\TTAab}{\put(0,90){\line(0,-1){30}}}
\newcommand{\TTAac}{\put(0,90){\line(3,-1){90}}}

\newcommand{\TTBaa}{\put(-90,60){\line(-1,-1){30}}}

\newcommand{\TTBac}{\put(-90,60){\line(1,-1){30}}}
\newcommand{\TTBbd}{\put(0,60){\line(-1,-1){30}}}

\newcommand{\TTBbf}{\put(0,60){\line(1,-1){30}}}
\newcommand{\TTBcg}{\put(90,60){\line(-1,-1){30}}}

\newcommand{\TTBci}{\put(90,60){\line(1,-1){30}}}

\newcommand{\TTCaa}{\put(-120,30){\line(-1,-3){10}}}

\newcommand{\TTCac}{\put(-120,30){\line(1,-3){10}}}

\newcommand{\TTCcg}{\put(-60,30){\line(-1,-3){10}}}

\newcommand{\TTCci}{\put(-60,30){\line(1,-3){10}}}
\newcommand{\TTCdj}{\put(-30,30){\line(-1,-3){10}}}

\newcommand{\TTCdl}{\put(-30,30){\line(1,-3){10}}}

\newcommand{\TTCfp}{\put(30,30){\line(-1,-3){10}}}

\newcommand{\BTAaL}[2]{\BTAa \put(-10,70){\makebox(20,20)[#1]{#2}}}

\newcommand{\BTBaL}[2]{\BTBa \put(-90,50){\makebox(20,20)[#1]{#2}}}
\newcommand{\BTBbL}[2]{\BTBb \put(70,50){\makebox(20,20)[#1]{#2}}}

\newcommand{\BTEaL}[2]{\BTEa \put(-160,-10){\makebox(20,20)[#1]{#2}}}

\newcommand{\BTEiL}[2]{\BTEi \put(0,-10){\makebox(20,20)[#1]{#2}}}

\newcommand{\TTAaL}[2]{\TTAa \put(-10,80){\makebox(20,20)[#1]{#2}}}

\newcommand{\TTBaL}[2]{\TTBa \put(-100,50){\makebox(20,20)[#1]{#2}}}
\newcommand{\TTBbL}[2]{\TTBb \put(-10,50){\makebox(20,20)[#1]{#2}}}
\newcommand{\TTBcL}[2]{\TTBc \put(80,50){\makebox(20,20)[#1]{#2}}}

\newcommand{\TTCaL}[2]{\TTCa \put(-130,20){\makebox(20,20)[#1]{#2}}}

\newcommand{\TTCcL}[2]{\TTCc \put(-70,20){\makebox(20,20)[#1]{#2}}}
\newcommand{\TTCdL}[2]{\TTCd \put(-40,20){\makebox(20,20)[#1]{#2}}}

\newcommand{\TTCfL}[2]{\TTCf \put(20,20){\makebox(20,20)[#1]{#2}}}
\newcommand{\TTCgL}[2]{\TTCg \put(50,20){\makebox(20,20)[#1]{#2}}}

\newcommand{\TTCiL}[2]{\TTCi \put(110,20){\makebox(20,20)[#1]{#2}}}

\newcommand{\TTDaL}[2]{\TTDa \put(-140,-10){\makebox(20,20)[#1]{#2}}}

\newcommand{\TTDcL}[2]{\TTDc \put(-120,-10){\makebox(20,20)[#1]{#2}}}

\newcommand{\TTDgL}[2]{\TTDg \put(-80,-10){\makebox(20,20)[#1]{#2}}}

\newcommand{\TTDiL}[2]{\TTDi \put(-60,-10){\makebox(20,20)[#1]{#2}}}
\newcommand{\TTDjL}[2]{\TTDj \put(-50,-10){\makebox(20,20)[#1]{#2}}}

\newcommand{\TTDlL}[2]{\TTDl \put(-30,-10){\makebox(20,20)[#1]{#2}}}

\newcommand{\TTDpL}[2]{\TTDp \put(10,-10){\makebox(20,20)[#1]{#2}}}

\begin{document}
\pagestyle{empty}
\title{The Wiener Polynomial of a Graph}
\author{Bruce E. Sagan \thanks{Supported in part by  grant
NSC 82-0501-1-001-A2
from the National
Science Council, Taiwan, R.O.C.} \\
\small{Department of Mathematics, Michigan State University,
East Lansing, MI 48824-1027, U.S.A.}\\[.1in]
Yeong-Nan Yeh \thanks{Supported in part by  grant 
NSC-83-0208-M001-68
from the National
Science Council, Taiwan, R.O.C.} \\
\small{Institute of Mathematics,
Academia Sinica,
Nankang, Taipei, Taiwan 11529, R.O.C.}\\[.1in]
and\\[.1in]
Ping Zhang\\
\small{Department of Mathematics, University of Texas,
El Paso, TX 79968-0514, U.S.A.}
}

\date{\today \\[1in]
	\begin{flushleft}
	\end{flushleft}
       }
\maketitle

\begin{flushleft} Proposed running head: \end{flushleft}
	\begin{center} 
The Wiener polynomial
	\end{center}

Send proofs to:

\begin{center}
\begin{tabular}{l}
Bruce E. Sagan\\
Department of Mathematics\\
Michigan State University\\ 
East Lansing, MI 48824-1027\\
U.S.A.\\
FAX: 517-432-1562\\
e-mail: sagan@math.msu.edu
\end{tabular}
\end{center}

	\begin{abstract}
The Wiener index is a graphical invariant that has found extensive
application in chemistry.  We define a generating function, which we
call the Wiener polynomial, whose derivative is a $q$-analog of the
Wiener index.  We study some of the elementary properties of this
polynomial and compute it for some common graphs.  We then find a
formula for the Wiener polynomial of a dendrimer, a certain highly
regular tree of interest to chemists, and show that it is
unimodal.  Finally, we point out a connection with the 
Poincar\'e polynomial of a finite Coxeter group.
	\end{abstract}
\pagestyle{plain}

\section{Introduction and elementary properties}		\label{iep}

Let $d(u,v)$ denote the distance between vertices $u$ and $v$ in a
graph $G$.  Throughout this paper we will assume that $G$ is connected.
The {\it Wiener index} of $G$ is defined as
	$$W(G)=\sum_{\{u,v\}}d(u,v)$$
where the sum is over all unordered pairs $\{u,v\}$ of distinct vertices in $G$.
The Wiener index was first proposed by Harold Wiener~\cite{wie:sdp}
as an aid to determining the boiling point of paraffin.  Since then,
the index has been shown to correlate with a host of other properties
of molecules (viewed as graphs).  For more information about the
Wiener index in chemistry and mathematics see~\cite{gp:mco}
and~\cite{bh:dg}, respectively.

We wish to define and study a related generating function.  If $q$ is
a parameter, then the {\it Wiener polynomial} of $G$ is
	\beq							\label{WGq}
	W(G;q)=\sum_{\{u,v\}}q^{d(u,v)}
	\eeq
where the sum is taken over the same set of pairs as before.
It is easy to see that the derivative of $W(G;q)$ is a $q$-analog of
$W(G)$ (see Theorem~\ref{pro}, number~\ref{qan}).  In the rest of this
section we will 
derive some basic properties of $W(G;q)$ and find its value when $G$
specializes to a number of simple graphs.  In Section~\ref{wpd} we
will compute the Wiener polynomial of a dendrimer $D_{n,d}$, a certain type of
highly regular tree that models various chemical molecules.  This
permits us to rederive results of Gutman and his
coauthors~\cite{gylc:wnd}.   We then use this formula to show that the
coefficients of $W(D_{n,d};q)$ are unimodal.  Finally we end with a
section on comments and open questions.  In particular, we point out
the connection with the  Poincar\'e polynomial of a Coxeter
group.

In what follows, any terms that are not defined will be found
described in the text of Chartrand and Lesniak~\cite{cl:gd}.  We will
use $|S|$ to denote the cardinality of a set $S$.  Also, if $f(q)$ is
a polynomial in $q$ then $\deg f(q)$ is its degree and $[q^i] f(q)$ is
the coefficient of $q^i$.  The next theorem summarizes some 
of the properties of $W(G;q)$.  Its proof follows easily from the
definitions and so is omitted.

\bth								\label{pro}
The Wiener polynomial satisfies the following conditions.
\ben
\item $\deg W(G;q)$ equals the diameter of $G$.
\item $[q^0] W(G;q)=0$.
\item $[q^1] W(G;q)=|E(G)|$ where $E(G)$ is the edge set of $G$.
\item $W(G;1)={|V(G)|\choose 2}$ where $V(G)$ is the vertex set of
			$G$. \label{WG1} 
\item $W'(G;1)=W(G).$\ \qed					\label{qan}
\een
\eth

We will next find the Wiener polynomial of some specific graphs.  
We let $K_n, P_n, C_n$ and $W_n$ denote the complete graph, path,
cycle and wheel on $n$ vertices, respectively.  Also let $Q_n$ be the
cube of dimension $n$ and $K_{m,n}$ be the complete bipartite graph on
parts of size $m$ and $n$.   Finally, $P$ denotes the Petersen graph.
Determining the Wiener polynomials of these graphs is a matter of
simple counting, so the proof of the next result is also omitted.
In the statement of the theorem we will use the standard $q$-analog of
$n$ which is $[n]=1+q+\cdots+q^{n-1}$.

\bth								\label{spe}
We have the following specific Wiener polynomials.
\ben
\item $W(K_n;q)={n\choose2}q$.
\item $W(K_{m,n};q)=mnq+\left[{m\choose2}+{n\choose2}\right]q^2$.
\item $W(W_n;q)=(2n-2)q+\frac{(n-1)(n-4)}{2}q^2$.
\item $W(P;q)=15q+30q^2$.
\item $W(P_n;q)=(n-1)q+(n-2)q^2+\cdots+q^{n-1}=\frac{q}{1-q}(n-[n])$.
					\label{WPn}
\item $W(C_{2n};q)=(2n)(q+q^2+\cdots+q^{n-1})+nq^n=2n([n]-1)+nq^n$.
\item $W(C_{2n+1};q)=(2n+1)(q+q^2+\cdots+q^n)=(2n+1)([n]-1)$.
\item $W(Q_n;q)=2^{n-1}[(1+q)^n-1]$.\ \qed
\een
\eth

Combining the previous theorem with number~\ref{qan} of
Theorem~\ref{pro}, we obtain the well-known Wiener indices of these
graphs.
\bth								\label{spe2}
We have the following specific Wiener indices.
\ben
\item $W(K_n)={n\choose2}$.
\item $W(K_{m,n})=(m+n)^2-mn-m-n$.
\item $W(W_n)=(n-1)(n-2)$.
\item $W(P)=75$.
\item $W(P_n)={n+1\choose3}$.
\item $W(C_{2n})=(2n)^3/8$.
\item $W(C_{2n+1})=(2n+2)(2n+1)(2n)/8$.
\item $W(Q_n)=n2^{2n-2}$.\ \qed
\een
\eth

It would be interesting to see what various graph
operations~\cite{har:gt} do to the 
Wiener polynomial.  Given graphs $G_1=(V_1,E_1)$ and $G_2=(V_2,E_2)$
with $|V_i|=n_i$ and $|E_i|=k_i$ for $i=1,2$ we define six new graphs
formed from $G_1,G_2$.  
\ben
\item {\it Join:} The graph $G_1+G_2$ has $V(G_1+G_2)=V_1\cup V_2$  and
$$E(G_1\cup G_2)=E_1\cup E_2\cup\{uv\ :\ u\in V_i, v\in V_2\}.$$
\een
In the other five cases the vertex set is always $V_1\times V_2$.
\ben
\item[2.] {\it Cartesian product:} The graph $G_1\times G_2$ has edge
set
$$\{(u_1,u_2)(v_1,v_2)\ :\ 
\mbox{$u_1v_1\in E_1$ and $u_2=v_2$ or $u_2v_2\in E_2$ and $u_1=v_1$}\},$$
\item[3.] {\it Composition:} The graph $G_1[G_2]$ has edge set
$$\{(u_1,u_2)(v_1,v_2)\ :\ 
\mbox{$u_1v_1\in E_1$ or $u_2v_2\in E_2$ and $u_1=v_1$}\},$$
\item[4.] {\it Disjunction:} The graph $G_1\jn G_2$ has edge set
$$\{(u_1,u_2)(v_1,v_2)\ :\ 
\mbox{$u_1v_1\in E_1$ or $u_2v_2\in E_2$ or both}\},$$
\item[5.] {\it Symmetric difference:} The graph $G_1\oplus G_2$ has edge set
$$\{(u_1,u_2)(v_1,v_2)\ :\ 
\mbox{$u_1v_1\in E_1$ or $u_2v_2\in E_2$ but not both}\},$$
\item[6.] {\it Tensor product:} The graph $G_1\times G_2$ has edge set
$$\{(u_1,u_2)(v_1,v_2)\ :\ 
\mbox{$u_1v_1\in E_1$ and $u_2v_2\in E_2$}\}.$$
\een
Taking a suggestion of Andreas Blass, it is
sometimes more natural to express our results in terms of the 
{\it ordered Wiener polynomial} defined by
	$$\Wb(G;q)=\sum_{(u,v)} q^{d(u,v)}$$
where the sum is now over all ordered pairs $(u,v)$ of vertices, including
those where $u=v$.  Thus
	$$\Wb(G;q)=2W(G;q)+|V(G)|.$$
Also it will be convenient to have a variable for the non-edges in
$G_i$ so let $\kb_i={n_i\choose2}-k_i$ for $i=1,2$.

\bpr							\label{G1G2}
Suppose $G_1$ and $G_2$ are connected and nontrivial (not equal to $K_1$).
Then with the preceding notation 
\ben
\item $W(G_1+G_2;q)=
(k_1+k_2+n_1n_2)q+(\kb_1+\kb_2)q^2$,
\item $\Wb(G_1\times G_2;q)=\Wb(G_1;q)\Wb(G_2;q)$,
\item $W(G_1[G_2];q)=n_1(k_2q+\kb_2q^2)+n_2^2W(G_1;q)$,
\item $W(G_1\jn G_2;q)=
(n_1^2k_2+n_2^2k_1-2k_1k_2)q+(n_1\kb_2+n_2\kb_1+2\kb_1\kb_2)q^2$,
\item $W(G_1\oplus G_2;q)=(n_1k_2+n_2k_1+2k_1\kb_2+2k_2\kb_1)q+
(n_1\kb_2+n_2\kb_1+2k_1k_2+2\kb_1\kb_2)q^2$.
\een
\epr
\pf\
In each part of this proof let $d_1,d_2$ and $d$ denote the distance
functions in $G_1,G_2$ and the graph formed from $G_1$ and $G_2$, respectively.

 1.  In $G_1+G_2$ all pairs of vertices are either at distance one or two.
If $d(u,v)=1$ then either $uv\in E_1$ or
$uv\in E_2$ or $u\in E_1, v\in E_2$.  This gives the linear
coefficient in $W(G_1+G_2;q)$.  The other term is gotten by counting
the remaining vertex pairs.

2.  A geodesic for an ordered pair $((u_1,u_2),(v_1,v_2))$ is obtained by
following a geodesic in $G_1$ for $(u_1u_2,v_1u_2)$ and then one
in $G_2$ for $(v_1u_2,v_1v_2)$.  The stated formula for $\Wb(G_1\times
G_2;q)$ follows.

3. First consider pairs $\{(u_1,u_2),(v_1,v_2)\}$ with $u_1=v_1$.  Then
$$
d((u_1,u_2),(u_1,v_2))=\case{1}{if $u_2v_2\in E_2$}{2}{if $u_2v_2\not\in E_2$}
$$
with a geodesic in the second case being $(u_1,u_2),(w_1,v_2),(u_1,v_2)$
where $w_1$ is any vertex adjacent to $u_1$ in $G_1$.
These vertex pairs contribute the first two terms in the sum for
$W(G_1[G_2];q)$.  If $u_1\neq v_1$ then
$d((u_1,u_2),(v_1,v_2))=d_1(u_1,v_1)$ since a $u_1$ to $v_1$ geodesic in
$G_1$ gives rise to one in $G_1[G_2]$ by adding a second component equal to
$u_2$ for the first vertex and to $v_2$ for all other vertices of the geodesic.

4.  This is similar to the previous proof where
$$
d((u_1,u_2),(v_1,v_2))=
\case{1}{if $u_1v_1\in E_1$ or $u_2v_2\in E_2$}{2}{else}
$$
with a geodesic in the second case being $(u_1,u_2),(w_1,w_2),(u_1,v_2)$
where $u_1w_1\in E_1$ and $w_2v_2\in E_2$.

5.  This is again similar to the previous two proofs with
$$
d((u_1,u_2),(v_1,v_2))=
\case{1}{if exactly one of $u_1v_1\in E_1$ or $u_2v_2\in E_2$}{2}{else.}
$$
In the second case, how to choose the middle vertex of the geodesic
depends on whether the neighborhoods of $u_i$ and $v_i$ are the same or
not, $i=1,2$,  as well as on whether both $u_1v_1\in E_1$ and $u_2v_2\in E_2$ or
neither. \hfill\Qed

As a corollary, we can rederive some of the $W(G;q)$
from Theorem~\ref{spe} as well as the ordered Wiener polynomial of the
grid $P_m\times P_n$.
\bco
We have the following specific Wiener polynomials
\ben
\item $W(K_{m,n};q)=mnq+\left[{m\choose2}+{n\choose2}\right]q^2$.
\item $W(W_n;q)=(2n-2)q+\frac{(n-1)(n-4)}{2}q^2$.
\item $\Wb(Q_n;q)=2^n(1+q)^n$.
\item $\Wb(P_m\times P_n)=\frac{1}{(1-q)^2}((1+q)m-2q[m])((1+q)n-2q[n])$.
\een
\eco
\pf\ For the first two polynomials use part 1 of  Theorem~\ref{G1G2}
as well as the fact that $W_n=C_{n-1}+K_1$ and $K_{m,n}=\Kb_m+\Kb_n$
where $\Kb_i$ is the completely disconnected graph on $i$ vertices.
For the last two polynomials use part 2 of the same theorem along with
the $n$-fold product $Q_n=K_2\times\cdots\times K_2$.  Note that 
in addition one needs $\Wb(K_2;q)=2+2q$ and
$\Wb(P_n;q)=2W(P_n;q)+n=\frac{1}{1-q}((1+q)n-2q[n])$. \hfill\Qed

Now we will find the Wiener polynomial  of a dendrimer, which will take
considerably more work.

\section{The Wiener polynomial of a dendrimer}			\label{wpd}

The {\it $d$-ary dendrimer on $n$ nodes}, $D_{n,d}$, is defined
inductively as follows.  The tree $D_{1,d}$ consists of a single node
labeled 1.  The tree $D_{n,d}$ has vertex set $\{1,2,\ldots,n\}$.  It 
is obtained  by attaching a leaf $n$ to the smallest numbered
node of $D_{n-1,d}$ which has degree  $\leq d$.  It is
convenient to consider $D_{n,d}$ as if it were rooted
at vertex number 1 with the nodes at each level ordered left to right
in increasing order of their numbering.  Thus in a typical tree 
the root has $d+1$ children while every other internal vertex 
(possibly with one exception) has $d$.  The dendrimer $D_{17,2}$ is
pictured in Figure~\ref{dT}

\thicklines
\setlength{\unitlength}{1.3pt}
\bfi
\bpi(180,130)(-130,-30)
\TTAaL{t}{\scz 1}
\TTBaL{l}{\mbc{\btab{c}\scz2\\ \scz(0,1)\st{18pt}\etab}}
\TTBbL{l}{\mbc{\btab{c}\scz3\\ \scz(1,0)\st{18pt}\etab}}
\TTBcL{r}{\mbc{\btab{c}\scz4\\ \st{18pt}\scz(1,1)\etab}}
\TTCaL{l}{\mbc{\btab{c}\scz5\\ \scz(0,1,0)\st{10pt}\etab}}
\TTCcL{l}{\mbc{\btab{c}\scz6\\ \scz(0,1,1)\st{10pt}\etab}}
\TTCdL{r}{\mbc{\btab{c}\scz7\\ \st{10pt}\scz(1,0,0)\etab}}
\TTCfL{r}{\mbc{\btab{c}\scz8\\ \scz(1,0,1)\etab}}
\TTCgL{r}{\mbc{\btab{c}\scz9\\ \scz(1,1,0)\etab}}
\TTCiL{r}{\mbc{\btab{c}\scz10\\ \scz(1,1,1)\etab}}
\TTDaL{l}{\mbc{\btab{c}\scz11\\ \scz(0,1,0,0)\etab}}
\TTDcL{b}{\mbc{\btab{c}\scz12\\ \scz(0,1,0,1)\etab}}
\TTDgL{l}{\mbc{\btab{c}\scz13\\ \scz(0,1,1,0)\etab}}
\TTDiL{bl}{\mbc{\btab{r}\scz14\\ \scz(0,1,1,1)\etab}}
\TTDjL{br}{\mbc{\btab{l}\scz15\\ \scz(1,0,0,0)\etab}}
\TTDlL{r}{\mbc{\btab{c}\scz16\\ \scz(1,0,0,1)\etab}}
\TTDpL{b}{\mbc{\btab{c}\scz17\\ \scz(1,0,1,0)\etab}}
\TTAaa
\TTAab
\TTAac
\TTBaa
\TTBac
\TTBbd
\TTBbf
\TTBcg
\TTBci
\TTCaa
\TTCac
\TTCcg
\TTCci
\TTCdj
\TTCdl
\TTCfp
\epi
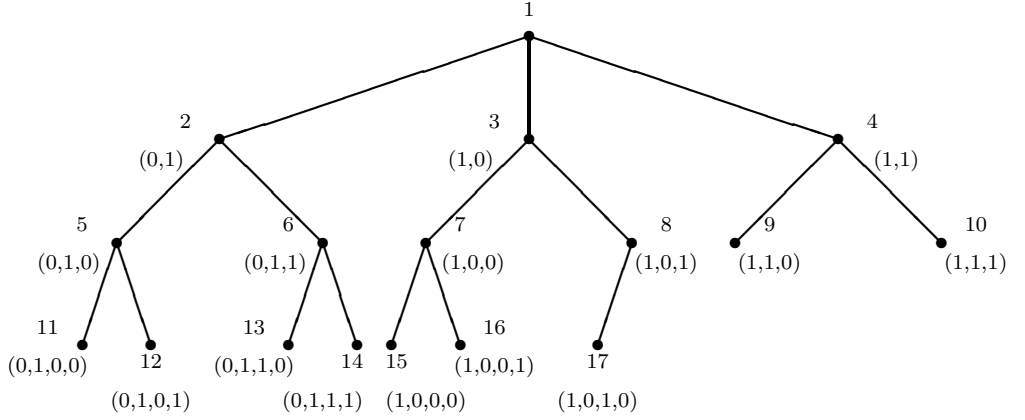
\capt{The dendrimer $D_{17,2}$}				\label{dT}
\efi

Define $n_k$ (respectively, $m_k$) to be the number of vertices in the
$d$-ary dendrimer  with exactly one descendant of vertex $2$
(respectively, of vertex $3$) at level $k+1$.  Thus
	$$n_k=2+(d+1)\frac{d^k-1}{d-1}\hs{.1in} \mbox{and}\hs{.1in} 
	m_k=3+(2d)\frac{d^k-1}{d-1}.$$
The tree in Figure~\ref{dT} has $n_0=2,n_1=5,n_2=11$ and
$m_0=3,m_1=7,m_2=15$. 

To describe $W(D_{n,d};q)$ we will also need to give each vertex $m>1$ a
label $\la(m)$ in addition to its number.  Specifically, if 
$n_k\leq m<n_{k+1}$ then
	$$\la(m)=(l_{k+1},l_{k},\cdots,l_0)\hs{.1in} \mbox{where}\hs{.1in}
	\sum_{i=0}^{k+1}l_id^i = n-n_k+(d-1)d^k,\
	0\leq l_i< d\ \forall i$$
so that the $l_i$ are the digits in the base $d$ expansion of
$n-n_k+(d-1)d^k$ (possibly with a leading zero).  Thus all the vertices at
level $k+1$ have labels 
which are consecutive base $d$ from $(0,d-1,\overbrace{0,\ldots,0}^{k})$
to $(1,\overbrace{d-1,\ldots,d-1}^{k+1})$.  This implies that if $m$
has label $\la(m)=(l_{k+1},l_{k},\cdots,l_0)$ then $m$'s children
are labeled left to right with $(l_{k+1},l_{k},\cdots,l_0,0)$ to
$(l_{k+1},l_{k},\cdots,l_0,d-1)$.
To illustrate, the labels of the vertices of $D_{17,2}$ are also given
in Figure~\ref{dT}.

To find $W(D_{n,d};q)$, we first 
consider the corresponding difference polynomials,
	$$\De W(D_{n,d};q)=W(D_{n,d};q)-W(D_{n-1,d};q).$$
So if $\De(W(D_{n,d};q))=\sum_i c_i q^i$ then $c_i$ is the number of
vertices in $D_{n,d}$ that are at distance $i$ from $n$.

\ble								\label{DeW}
Suppose $\la(n)=(l_{k+1},l_{k},\cdots,l_0)$.  Then
	$$\De W(D_{n,d};q)=
	\sum_{i=0}^k d^i q^{2i+1} + \sum_{i=0}^{k'} d^i(l_i+1)q^{2i+2}$$
where $k'=k-1$ or $k$ for $n_k\leq n< m_k$ or
$m_k\leq n< n_{k+1}$, respectively.
\ele
\pf\  We will do the case  $n_k\leq n< m_k$, the other being similar.
Suppose first that $n=n_k$ so that $l_0=\ldots=l_{k-1}=0$.  Thus we
wish to show
	$$\De W(D_{n_k,d};q)=q+q^2+dq^3+dq^4+\cdots+d^kq^{2k+1}.$$
Every vertex at distance $i-1$ from $n_{k-1}$ in $D_{n_{k-1},d}$ is at
distance $i$ from $n_k$ in $D_{n_k,d}$.  When $i=2j$ is even, this
accounts for all the vertices at distance $2j$ from $n_k$ and so
	$$[q^{2j}]\De W(D_{n_k,d};q)=
	[q^{2j-1}]\De W(D_{n_{k-1},d};q)=d^{j-1}$$
by induction.  When $i=2j+1>1$ then any leaves of $D_{n_k,d}$ which
are descendants of $n_{k-j-1}$ but not of $n_{k-j}$ are also at
distance $2j+1$ from $n_k$.  (If $j=k$ then let $n_{-1}=1$.)  So
by induction
	\bea
	[q^{2j+1}]\De W(D_{n_k,d};q)&=&
		[q^{2j}]\De W(D_{n_{k-1},d};q)+
			\mbox{\#(leaves at distance $2i+1$)}\\
	&=&\case{d^{j-1}+d^{j-1}}{if $j<k$}{0+d^k}{if $j=k$}\\
	&=&d^j\ \mbox{for all $j$.}
	\eea
This completes the proof when $n=n_k$.

Now suppose that $n_k<n<m_k$.  
We will construct a function $f$ that sets up a bijection between
vertices at distance $i$ from 
vertex $n$ in $D_{n,d}$ and those at distance $i$ from vertex $n'$
in $D_{n',d}$ where $n_k\leq n'<n$ and then use induction.
The bijection will hold for all $i$, $0\leq i\leq 2k+1$, except
for $i=2j$ where $j$ will be determined shortly.

All vertices satisfying $n_k<n<m_k$ are
descendents of vertex 2.  However, since $n\neq n_k$ the unique
$2$ to $n$ path contains an edge $uv$ where $v$ is not the leftmost child of
$u$.  Let $u$ be the lowest such vertex and define
	\bea
	j&=&\mbox{distance from $u$ to $n$,}\\
	v'&=&\mbox{child of $u$ just to the left of $v$,}\\
	n'&=&\mbox{leftmost lowest descendant of $v'$.}
	\eea
For a schematic sketch of this situation for $d=2$, see Figure~\ref{uvv}.

\bfi
\bpi(300,130)(-150,-30)
\BTAaL{t}{$u$}
\BTBaL{t}{$v'$}
\BTBbL{t}{$v$}
\BTDa
\BTDb
\BTDd
\BTDe
\BTDh
\BTEaL{b}{$n'$}
\BTEb
\BTEc
\BTEd
\BTEg
\BTEh
\BTEiL{b}{$n$}
\BTAaa
\BTAab
\BTBaa
\BTBab
\BTBbc
\BTBbd
\BTCaa
\BTCbd
\BTCce
\BTCdh
\BTDaa
\BTDab
\BTDbc
\BTDbd
\BTDdg
\BTDdh
\BTDei
\put(-140,20){\line(1,0){120}}
\put(20,20){\line(1,0){120}}
\put(-130,-20){\line(0,1){10}}
\put(-10,-20){\line(0,1){10}}
\put(-130,-20){\line(1,0){120}}
\put(-80,40){\makebox(0,0){$D_{n',d}(v')$}}
\put(80,40){\makebox(0,0){$D_{n,d}(v)$}}
\put(-60,10){\makebox(0,0){$\cdots$}}
\put(-70,-30){\makebox(0,0){$\in D_{n,d}\setminus D_{n',d}$}}
\epi
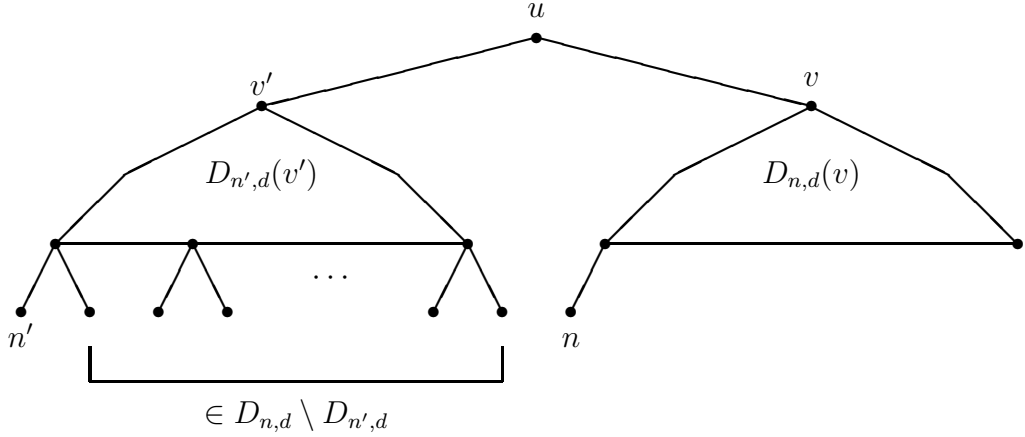
\capt{Location of $u,v,v',n$ and $n'$ when $d=2$}				\label{uvv}
\efi

It follows from the definitions and the way in which children of a vertex
are labeled that
	\beq							\label{lan}
	\la(n)=(l_{k+1},\ldots,l_0)=(l_{k+1},\ldots,l_j,l,0,\ldots,0)
	\eeq
and
	\beq							\label{lan'}
	\la(n')=(l'_{k+1},\ldots,l'_0)=(l_{k+1},\ldots,l_j,l-1,0,\ldots,0)
	\eeq	
where $l>0$.  So $l_i=l'_i$ if $i\neq j-1$ and $l_{j-1}$ only enters
into the coefficient of $q^{2j}$.  Construction of the bijection $f$
will establish that
$[q^i]\De W(D_{n,d};q)=[q^i]\De W(D_{n',d};q)$ for $i\neq 2j$.
Thus induction combined with the equations for $\la(n)$ and $\la(n')$
will complete the proof in this case.
	
Construct $f$ as follows.  Given any ordered tree $T$ and one of its vertices
$v$, then we let $T(v)$ denote the subtree of $T$ consisting of $v$
and all its descendants.  Note that there is a unique isomorphism of
ordered trees $g:\ D_{n,d}(v)\ra D_{n',d}(v')$.  Also note that
all vertices of $D_{n,d}\setminus D_{n',d}$ (other than $n$
itself) are leaves at distance $2j$ from $n$.  So if $w$ is any vertex
of $D_{n',d}\cup n$ then let
	$$f(w)=\left\{\barr{ll}
		g(w)		&\mbox{if $w\in D_{n,d}(v)$}\\
		g^{-1}(w)	&\mbox{if $w\in D_{n',d}(v')$}\\
		w		&\mbox{else.}
		\earr\right.$$
Thus $d(n,w)=d(n',f(w))$:  For $w\in D_{n,d}(v)\cup D_{n',d}(v')$
this follows because $g$ is an isomorphism.  For any other $w$, the
unique $n$ to $w$ and $n'$ to $w$ paths both go through $u$, so
$d(n,w)=d(n',w)$.  The function $f$ is also clearly bijective, so we
are done when $i\neq 2j$.

To complete the proof, note that $f$ restricts to an injection from
the vertices at distance $2j$ from $n'$ in $D_{n',d}$
into those at distance $2j$ from $n$ in $D_{n,d}$.  The only
remaining $w$ with $d(w,n)=2j$ are the leaves of $D_{n,d}(v')$ which
are $d^{j-1}$ in number.  So by induction and equations~(\ref{lan})
and~(\ref{lan'})
	$$[q^{2j}]\De W(D_{n,d};q)=[q^{2j}]\De W(D_{n',d};q)+d^{j-1}
		=d^{j-1}(l+1)$$
as desired. \hfill\Qed

We are now in a position to compute $W(D_{n,d};q)$.  In the
following theorem $\lfloor\cdot\rfloor$ denotes the floor (round down)
function.  
\bth								\label{WTn}
Suppose $\la(n)=(l_{k+1},l_{k},\cdots,l_0)$ and define
	$$\la_i(n)=l_0 + l_1 d + \cdots l_{i-1} d^{i-1}.$$  
If $n_k\leq n<m_k$ then
	\bea
	W(D_{n,d};q)&=&
	\sum_{i=0}^k d^i(n-n_i+1)q^{2i+1} + \\
	&&\hs{-25pt}
	\sum_{i=0}^{k'}\left(d^{2i}\flf{n-m_i}{d^{i+1}}{d+1\choose2}+
	d^{2i}{l_i+1\choose2}+d^i(l_i+1)(\la_i(n)+1)\right)q^{2i+2}.	
	\eea
where $k'=k-1$ or $k$ for $n_k\leq n< m_k$ or
$m_k\leq n< n_{k+1}$, respectively.
\eth
\pf\  Since $W(D_{1,d};q)=0$, we have 
$W(D_{n,d};q)=\sum_{m=2}^n\De	W(D_{m,d};q)$.  Thus the theorem
will follow from Lemma~\ref{DeW} and summation.  Summing the
coefficients of $q^{2i+1}$ in the various $\De W(D_{m,d};q)$
is easy since the nonzero ones are all equal to $d^i$.
The nonzero  $[q^{2i+2}]\De W(D_{m,d};q)$ are periodic with period
$d^{i+1}$.  In fact this periodic sequence is 
	$$d^i,\cdots,d^i,2d^i,\cdots,2d^i,\cdots,d^{i+1},\cdots,d^{i+1}$$
where each block of equal values is of size $d^i$.  Thus the sum of a
complete period is $d^{2i}{d+1\choose2}$ and the term containing this
expression in $[q^{2i+2}] W(D_{n,d};q)$ comes from summing all
complete periods prior to the (perhaps partial) period containing 
$c=[q^{2i+2}] \De W(D_{n,d};q)$. The term $d^{2i}{l_i+1\choose2}$
comes from summing all complete blocks of the last period prior to the
(perhaps partial) block containing $c$.  Finally, 
$d^i(l_i+1)(\la_i(n)+1)$ is the contribution of the block containing
$c$. \hfill\Qed

We will now rederive the Wiener index of a complete dendrimer, as was
first done in~\cite{gylc:wnd}.  A {\it complete dendrimer} is $D_{n,d}$
where $n=n_k-1$, i.e., as an ordered tree it is complete to level $k$.
The reader should be warned that in~\cite{gylc:wnd} they only consider
complete dendrimers and index them with two parameters that are
different from ours.  We have
	$$\la(n_k-1)=(1,\overbrace{d-1,\ldots,d-1}^{k})$$
and
	$$\la_i(n_k-1)=d^i-1$$
for $i<k$.  Substituting these values into the formulas of
Theorem~\ref{WTn} gives
	\bea
	W(D_{n_k-1,d};q)&=&
	\sum_{i=0}^{k-1}d^i\cdot d^i(d+1)\frac{d^{k-i}-1}{d-1}q^{2i+1}\\
    &&+\sum_{i=0}^{k-1}\left[d^{2i}(d+1)\frac{d^{k-i-1}-1}{d-1}{d+1\choose2}+
		d^{2i}{d\choose2}+d^i\cdot d\cdot d^i\right]q^{2i+2}\\
    &=&\sum_{i=0}^{k-1}d^{2i}(d+1)\frac{d^{k-i}-1}{d-1}q^{2i+1}\\
    &&+\sum_{i=0}^{k-1}d^{2i}{d+1\choose2}
		\left[(d+1)\frac{d^{k-i-1}-1}{d-1}+1\right]q^{2i+2}.
	\eea
Taking the derivative of the previous equation and setting $t=1$ gives
the Wiener index according to Theorem~\ref{pro}, number~\ref{qan}.  To
evaluate the summations, use
	$$\sum_{i=0}^{k-1}(i+1)d^i=\frac{kd^{k+1}-(k+1)d^k+1}{(d-1)^2}$$ 
and its variants repeatedly to obtain
	\bea
	W(D_{n_k-1,d})&=&\frac{d+1}{d-1}
	\left[\frac{(2k-1)d^{2k+1}-(2k+1)d^{2k}+d^{k+1}+d^k}{(d-1)^2}\right.\\
        &&\hspace*{9pt}
	-\left.\frac{(2k-1)d^{2k+2}-(2k+1)d^{2k}+d^{2}+1}{(d^2-1)^2}\right]\\
	&&+\frac{d(d+1)^2}{d-1}
	\left[\frac{kd^{2k}-(k+1)d^{2k-1}+d^{k-1}}{(d-1)^2}\right.\\
        &&\hspace*{9pt}
	-\left.\frac{kd^{2k+2}-(k+1)d^{2k}+1}{(d^2-1)^2}\right]
	+d(d+1)\frac{kd^{2k+2}-(k+1)d^{2k}+1}{(d^2-1)^2}.
	\eea
Simplification of the above expression yields the following result
which is equivalent to equation (9) of~~\cite{gylc:wnd} after a change
of variables.

\bco
The Wiener index of a complete dendrimer is
	$$W(D_{n_k-1,d})=
  \frac{d^{2k}[kd^3+(k-2)d^2-(k+3)d-(k+1)]+2d^k(d+1)^2-(d+1)}{(d-1)^3}.$$
\eco

\section{Unimodality}						\label{u}

We say a sequence $(a_m)_{m\geq0}$ is {\it unimodal} if, for some
index $k$
	$$a_0\leq a_1\leq \ldots \leq a_k\geq a_{k+1}\geq a_{k+2}\geq\ldots.$$
Unimodal sequences appear in many areas of mathematics.  For a survey, see
Stanley's article~\cite{sta:lcu}.  We will show that the coefficients
of $W(D_{n,d})$ are
unimodal.  First, however, we will need a general result about
sequences and their differences.    The {\it difference sequence} 
of $(a_m)_{m\geq0}$ is $(\De a_m)_{m\geq1}$ where $\De a_m= a_m-a_{m-1}$.

\bpr								\label{am}
Suppose $(a_m)_{m\geq0}$ and $(b_m)_{m\geq0}$ are two sequences with
$\De a_m\leq \De b_m$ for all $m$ in some interval $I$ of integers.
If $a_M\leq b_M$ for some $M\in I$ then $a_m\leq b_m$ for all $m\geq M$
with $m\in I$.  On the other hand, if $a_M\geq b_M$ for some $M\in I$
then $a_m\geq b_m$ for all $m\leq M$ with $m\in I$.
\epr
\pf\  We will consider the case $a_M\leq b_M$, the other being
similar.  Using the given inequalities we have, for $m\geq M$,
	$$a_m = a_M + \sum_{k=M+1}^m \De a_k\leq
		b_M + \sum_{k=M+1}^m \De b_k = b_m$$
as desired. \hfill\Qed

We will need coefficients of $W(D_{n,d};q)$ to play the roles of
$a_M$ and $b_M$ in the previous Proposition.
\ble								\label{pk}
Let $p_k=m_k+2d^k-1$ and set $W(t)=W(D_{p_k,d};q)$.  Then
	$$[q^{2k}]W(t)\leq [q^{2k+1}]W(t) = [q^{2k+2}]W(t).$$
\ele
\pf\ Note that by the choice of $p_k$ we have
	\beq							\label{lap}
	\la(p_k)=(1,1,\overbrace{d-1,\ldots,d-1}^k)
	\eeq
and so
	\beq							\label{lai}
	\la_{i}(p_k)=d^i-1
	\eeq
for $i\leq k$.  Using Theorem~\ref{WTn} and the previous two equations we get
	$$[q^{2k}]W(t)=
	d^{2k-2}\flf{p_k-m_{k-1}}{d^{k}}{d+1\choose2}+
	d^{2k-2}{d\choose2}+d^{k-1}\cdot d\cdot d^{k-1}.$$
But $p_k-m_{k-1}=m_k-m_{k-1}+2d^k-1=4d^k-1$ so
	\beq							\label{t2k}
	[q^{2k}]W(t)=
	d^{2k-1}\left(3\frac{d+1}{2}+\frac{d-1}{2}+1\right)=2d^{2k-1}(d+1)
	\eeq

Using Theorem~\ref{WTn} again we obtain
	\beq							\label{t2k+}
	[q^{2k+1}]W(t)=d^k(p_k-n_k+1)=d^k(m_k-n_k+2d^k)=3d^{2k}.
	\eeq
Comparison of equations~(\ref{t2k}) and~(\ref{t2k+}) 
and the fact that $d\ge2$ yield the 
inequality in the statement of lemma.  

One last application of Theorem~\ref{WTn} together with
equations~(\ref{lap}) and~(\ref{lai}) gives
	$$[q^{2k+2}]W(t)=
	d^{2k}\flf{p_k-m_{k}}{d^{k+1}}{d+1\choose2}+
	d^{2k}{2\choose2}+d^{k}\cdot 2\cdot d^{k}=0+d^{2k}+2d^{2k}=3d^{2k}.$$
But this is the same as the value obtained in equation~(\ref{t2k+}), so
we are done with the proof. \hfill\Qed

We now put the various pieces together to get the promised theorem.
\bth
Let $W(t)=W(D_{n,d};q)$.  Then the coefficients of $W(t)$ are
unimodal.  Furthermore, for $n_k\leq n< p_k$ either $[q^{2k}]W(t)$ or
$[q^{2k+1}]W(t)$ is a maximum.  On the other hand, for $p_k\leq
n<n_{k+1}$ we have a maximum at $[q^{2k+2}]W(t)$ so that in this case
the sequence is increasing.
\eth
\pf\  We will only consider the case $n_k\leq n< p_k$.  The reader
can easily fill in the details in the other one.  Since 
$\deg W(t)\leq 2k+2$ it suffices to show that the following two
equations hold:	
	\bea
	[q^i]W(t)&\leq&[q^{i+1}]W(t)\ \mbox{for $i<2k$, and} \\  
	\left[q^{2k+1}\right]W(t)&\geq&[q^{2k+2}]W(t)
	\eea

For the first inequality, fix $i$ and consider the sequence whose
terms are  given by $a_m=[q^i]W(D_{m,d};q)$.  So in particular
$a_n=[q^i]W(t)$.  Similarly define $b_m=[q^{i+1}]W(D_{m,d};q)$ so
that $b_n=[q^{i+1}]W(t)$.  Lemma~\ref{DeW} shows that these two
sequences satisfy the supposition of Proposition~\ref{am} for 
$m\geq m_{\lfloor i/2 \rfloor}$.  By the lemma just proved,
$a_{p_{\lfloor i/2 \rfloor}}\leq b_{p_{\lfloor i/2 \rfloor}}$.  Since
$n\geq n_k>p_{\lfloor i/2 \rfloor}$, Proposition~\ref{am} applies to
show that $a_n\leq b_n$ as desired.

The second inequality is clearly true for $n_k\leq n< m_k$ since then 
$[q^{2k+2}]W(t)=0$.  If $m_k\leq n<p_k$, then consider the sequences
defined by $a_m=[q^{2k+1}]W(D_{m,d};q)$ and
$b_m=[q^{2k+2}]W(D_{m,d};q)$.   These satisfy the hypothesis of 
Proposition~\ref{am} for $m\geq m_k$.  By the previous lemma
$a_{p_k}=b_{p_k}$.  So, since $n<p_k$, we can apply
Proposition~\ref{am} to get $a_n\geq b_n$ which concludes the proof. \hfill\Qed

\section{Comments and open questions}

Since the Wiener polynomial is a new graphical invariant, there are
many questions one could ask about it.  We summarize some of them in
this section

(I) The reader has probably noticed that we did not provide a formula
for $W(G_1\times G_2;q)$ in terms of $W(G_1;q)$ and $W(G_2;q)$.  It
would be interesting to fill this gap in the list of Wiener
polynomials related to graph operations.

\medskip

(II)  A referee pointed out that the generating function for the
Wiener polynomial of the complete 
dendrimer has a nice form.  It can be obtained algebraically from the
equation for $W(D_{n_k-1,d};q)$ given at the end of
section~\ref{wpd}, but the referee asked for a combinatorial proof.
We give such a demonstration next.
\bpr  The generating function for $W(D_{n_k-1,d};q)$ is
$$
\sum_{k\ge1} W(D_{n_k-1,d};q) z^k=
z\frac{(d+1)q+{d+1\choose 2}q^2(1+z)}{(1-z)(1-dz)(1-d^2q^2z)}.
$$
\epr
\pf\  Let $T_k=D_{n_k-1,d}$ and $F(z)=\sum_{k\ge1} W(T_k;q) z^k$.
Consider $T_{k-1}$ as embedded in $T_k$ so that their roots coincide.
Then $(1-z)F(z)$ is the generating function for all $u$ to $v$ paths in
$T_k$ such that at least one of $u,v$ is a leaf, $u\neq v$.  Now there is a
$d^2$-to-1 mapping from paths $P$ of length $l$ in $T_k$ to paths $P'$ of
length $l-2$ in $T_{k-1}$ gotten by removing the two endpoints of $P$.
Furthermore one of the endpoints of $P$ is a leaf iff one of the
endpoints of $P'$ is a leaf.  So $(1-d^2q^2z)(1-z)F(z)$ is the
generating function for all paths in $T_k$ of length one or two with
at least one endpoint a leaf.  Now if $k>2$ then $T_k$ has exactly
$d$ times as many such paths as $T_{k-1}$.  So 
$(1-dz)(1-d^2q^2z)(1-z)F(z)$ is a polynomial of degree 2 which is easy
to compute directly, giving the numerator of the fraction in the
statement of the proposition.\hfill\qed

\medskip

(III) The Wiener polynomial refines the Wiener index since it gives
information about how many pairs of vertices are at a given distance
$i$, not just the sum of all distances.  One can also refine the
Wiener polynomial itself as follows.  Define the {\it Wiener
polynomial of a graph $G$ relative to a vertex $v$} by
	$$W_v(G;q)=\sum_{w\in V} q^{d(v,w)}.$$
The next result is immediate from the definitions and will be useful
latter.  

\bpr								\label{WvG}
We have the following relationship between the ordered and relative
Wiener polynomials
	$$\Wb(G;q)=\sum_{v} W_v(G;q).$$
Furthermore, if $G$ is vertex transitive then for any vertex $v\in G$
we have
	$$\Wb(G;q)=|V| W_v(G;q).\ \Qed$$
\epr

\medskip

(IV)  In certain cases, Wiener polynomial is closely related to a
polynomial that appears in the theory of Coxeter groups.  We will
follow the text of Humphreys~\cite{hum:rgc} in terms of definitions
and notation.  Let $(W,S)$ be a Coxeter system and let 
$T=\{wsw^{-1}\ :\ w\in W,s\in S\}$.  There should be no confusion
between the notation for a Coxeter group $W$ and the one for the
Wiener invariants since the latter is always followed by a
parenthesized expression.  The {\it absolute length} of 
$w\in W$ is the minimum number $k$ such that
	$$w=t_1 t_2 \cdots t_k\hs{10pt}\mbox{where $t_i\in T\ \forall i.$}$$
We write $\lh(w)=k$ in this case.  The absolute length function differs
from the ordinary length function in that the factors in the product
for $w$ are required to be in $T$ rather than in $S$.  Recently
Barcelo, Garsia and Goupil~\cite{bg:nbc} have found a beautiful
connection between absolute length and NBC bases.

Our interest is in the {\it  Poincar\'e polynomial} of $W$
which is defined by
	$$\Pi(W;q)=\sum_{w\in W} q^{\lh(w)}.$$
It is related to the extended Wiener polynomial as follows.  We will define a
graph $G_W$ associated with any Coxeter group.  The vertices of $G_W$
are the elements of $W$ and we connect $v$ and $w$ by an edge if
$v=wt$ for some $t\in T$.  This  graph is related to the
strong Bruhat ordering of $W$.  Note that $\lh(w)=d(1,w)$ where 1 is
the identity element of $W$ and distance is taken in $G_W$.  Combining
this observation with Proposition~\ref{WvG} and the easily proved fact
that $G_W$ is vertex transitive, we obtain the connection between
the two polynomials.
\bpr								\label{WGW}
If $W$ is a finite Coxeter group and $G_W$ is the corresponding graph then
	$$\Wb(G_W;q)=|W|\ \Pi(W;q).\ \qed$$
\epr

The next theorem is well known, see the book of Orlik
and Terao~\cite{ot:ah}.
\bth								\label{APW}
If $W$ is a finite Coxeter group then its Poincar\'e
polynomial factors as
	$$\Pi(W;q)=\prod_{e} (1+eq)$$
where the product is over all exponents $e$ of $W$.  In particular,
the roots of $\Pi(W;q)$ are all negative rational numbers.\ \hfill\Qed
\eth
The comment about the roots of $\Pi(W;q)$ relates to the concept of
unimodality as follows.  We say a sequence $(a_m)_{m\geq0}$ is {\it
log concave} if $a_m^2\geq a_{m-1}a_{m+1}$ for all $m\geq1$.  The
relationship between these three concepts is easy to prove
(see~\cite{sta:lcu}). 
\bpr								
Let $(a_m)_{0\le m\le M}$ be a  sequence of positive numbers
and let $f(q)=\sum_{m=0}^M a_mq^m$ be the corresponding polynomial. 
\ben
\item If $f(q)$ factors over the negative rationals then $(a_m)$ is
log concave.
\item If $(a_m)$ is log concave then it is also unimodal.\ \Qed
\een
\epr
This brings up three questions.
\ben
\item  For which graphs $G$ is the coefficient sequence of $W(G;q)$
unimodal?
\item  For which graphs $G$ is the coefficient sequence of $W(G;q)$
log concave?
\item  For which graphs $G$ does $\Wb(G;q)$ factor over the negative rationals.
\een
Note that both the graphs $G_W$ and $Q_n$ satisfy the last condition
which is the strongest of the three.

\medskip

(V)  There are two theorems that are useful for computing the Wiener
index of a tree that we have been unable to find analogs for in the
case of the Wiener polynomial.  The first is due to Wiener
himself~\cite{wie:sdp}.
\bth
Let $T$ be a tree and let $e$ be an edge in $E=E(T)$.  Let $n_1(e)$ and
$n_2(e)$ be the number of vertices in the two components of $T-e$.
Then
	$$W(T)=\sum_{e\in E} n_1(e)n_2(e).\ \Qed$$
\eth

The second is a result of Gutman~\cite{gut:nmc}.  It is useful when
the tree in question has few vertices of high degree.
\bth
Let $T$ be a tree and let $v$ be a vertex in $V=V(T)$ with 
$\deg v\geq3$.  Let
$n_1(v),\ldots,n_{\deg v}(v)$ be the number of vertices in each of the
components of $T-v$.  Then
	$$W(T)={|V|+1\choose3}-
  \sum_{v\in V}\sum_{1\leq i<j<k\leq\deg v} n_i(v)n_j(v)n_k(v).\ \Qed$$
\eth
Note that this theorem immediately gives the Wiener number of the path
$P_n$ as in Theorem~\ref{spe}, number~\ref{WPn}, as well as the fact
that $P_n$ has the largest Wiener number of any tree.  It would be
nice to get analogous results for Wiener polynomials, possibly be
comparing them coefficient-wise.

\medskip

(VI)  Because of the development of parallel architectures for
interconnection computer networks, there has recently been interest in
a generalization of the distance concept.  A {\it container},
$C(u,v)$, is a set of
vertex-disjoint paths between two vertices $u,v\in V(G)$, i.e., any
two paths in $C(u,v)$ only intersect at $u$ and $v$.  The {\it width},
$w=w(C(u,v))$, is the number of paths in the container while the {\it
length}, $l=l(C(u,v))$, is the length of the longest path in $C(u,v)$.
For fixed $w$, define the {\it $w$-distance between $u$ and $v$} as
	$$d_w(u,v)=\min_{C(u,v)} l(C(u,v)).$$
where the minimum is taken over all containers $C(u,v)$ of width $w$.
Note that when $w=1$ then $d_w(u,v)$ reduces to the usual distance
between $u$ and $v$.
For more information about these concepts and their relation to
networks, see the article of Hsu~\cite{hsu:cwl}.

Now we can define the $w$-Wiener polynomial by
	$$W_w(G;q)=\sum_{u,v} q^{d_w(u,v)}.$$
It would be interesting to compute this polynomial for various graphs
and study its properties, e.g., unimodality.  It would also be
interesting to see if this object yields any useful information in
chemistry, group theory, or computer science.

\begin{\bib}{99}

\bibitem{bg:nbc}  H. Barcelo and A. Goupil, Non broken circuits of
reflection groups and their factorization in $D_n$, 
{\it Israel J. Math.} {\bf 91} (1995), 285--306.

\bibitem{bh:dg} F. Buckley and F. Harary, ``Distance in Graphs,''
Addison-Wesley, Redwood, CA, 1990.

\bibitem{cl:gd} G. Chartrand and L. Lesniak, ``Graphs and Digraphs,'' 
second edition, Wadsworth \& Brooks/Cole, Monteray, CA, 1986.

\bibitem{gut:nmc}  I. Gutman, A new method for the calculation of the
Wiener number of acyclic molecular graphs, 
{\it Journal of Molecular Structure (Theochem)} {\bf 285} (1993), 137--142.

\bibitem{gp:mco} I. Gutman and O. Polansky, {\it Mathematical Concepts
in Organic Chemistry}, Springer-Verlag, Berlin, Germany, 1986.

\bibitem{gylc:wnd} I. Gutman, Y.-N. Yeh, S.-L. Lee J.-C. Chen, Wiener
numbers of dendrimers,
{\it Comm.\  Math.\ Chem.\ } {\bf 30} (1994), 103--115.

\bibitem{har:gt} F. Harary, ``Graph Theory,'' Addison-Wesley, Reading,
MA, 1971.

\bibitem{hsu:cwl} D. F. Hsu, On container width and length in graphs,
groups and networks, {\it IEICE Trans.  on Fundamentals of
Electronics, Communications and Computer Science}, 
{\bf Vol.\ E77-A, No.\  4} (1994), 668--680.

\bibitem{hum:rgc} J. E. Humphreys, ``Reflection Groups and Coxeter
Groups,'' Cambridge Studies in Advanced Mathematics, Cambridge
University Press, Cambridge, 1990.

\bibitem{ot:ah} P. Orlik and H. Terao, ``Arrangements of Hyperplanes,''
Grundlehren 300, Springer-Verlag, New York, NY, 1992.

\bibitem{sta:lcu} R. P. Stanley, Log-concave and unimodal 
sequences in algebra, combinatorics, and geometry, in
``Graph Theory and Its Applications: East and West,'' Ann. NY
Acad. Sci. {\bf 576} (1989), 500--535.

\bibitem{wie:sdp} H. Wiener, Structural determination of paraffin
boiling points, {\it J. Amer. Chem. Soc.} {\bf 69} (1947), 17--20.

\end{\bib}

\end{document}